\documentclass[a4,11pt]{article}
\usepackage{amsmath,amsthm,amssymb}
\newtheorem{Th}{Theorem}[subsection] 
\newtheorem{Lem}[Th]{Lemma} 
\newtheorem{Prop}[Th]{Proposition} 
\newtheorem{Cor}[Th]{Corollary} 
 
\newtheorem{Prob}[Th]{Problem}

\newtheorem{Ex}[Th]{Example}

\def\N{\mathbb{N}} 
\def\R{\mathbb{R}}

\def\E{{\cal E}} 
\def\cH{{\cal H}}

\def\cO{{\cal O}}

\def\U{{\cal U}}

\def\cN{{\cal N}} 
\def\cL{{\cal L}}

\def\half{\frac{1}{2}}

\def\DIV{\mbox{\rm div}} 
\def\supp{\mbox{\rm supp}}

\def\ov{\overline}

\def\1/2{\frac{1}{2}}

\def\srT{{\mbox{\tiny\rm T}}}

\makeatletter 
\long\def\@makefntext#1{\parindent 1em\noindent 
\@hangfrom{\hbox to 1.8em{\hss$^{\@thefnmark}$}}#1}
\makeatother

\begin{document}
\mbox{}
\vspace{-0.5cm}
\begin{center}
{\large\bf Shape derivative of potential energy\\
and energy release rate in fracture mechanics}\\
~\\
~\\
{Masato Kimura$^{*}$
and Isao Wakano$^{**}$}\\
~\\
{\small
$^*$ {\it Faculty of Mathematics, Kyushu University,
masato@math.kyushu-u.ac.jp}\\
$^{**}$ {\it Graduate School of Informatics, Kyoto University,
wakano@i.kyoto-u.ac.jp}
}
\end{center}

~\\

\begin{abstract} 
We study general mathematical framework for 
variation of potential energy with respect to 
domain deformation.
It enables rigorous derivation of the integral formulas for
the energy release rate in crack problems.
Applying a technique of
the shape sensitivity analysis, we formulate the 
shape derivative of potential energy as 
a variational problem with a parameter. Key tools of
our abstract theory are a new parameter variational principle  
and the classical implicit function theorem in Banach
spaces.
\end{abstract}

\subsection{Introduction}\label{intro} 
\setcounter{equation}{0}

Many variational problems related to the domain
deformation have been investigated in the
theory of shape derivatives or the shape sensitivity
analysis (\cite{Kim00}, \cite{S-Z92}).
According to the spread of the importance of shape derivatives
in various scientific fields, more development of their mathematical
foundation has been required.
The purpose of this paper is
to establish some abstract parameter variational formulas
and their application to the shape derivative of potential energies.
An important example is the energy release rate in crack problems,
which is known as one of the most fundamental quantities
in the theory of fracture mechanics.  

Scientific investigation to understand 
crack evolution process in elastic body
was originated by Griffith \cite{Grif20}
and has been studied from various viewpoints
in engineering, physics and mathematics since then.
Griffith's idea in the fracture mechanics is
even now the fundamental theory in 
modeling and analysis of the crack behavior.
We here make reference to only very few extended studies
from mathematical point of view,
Cherepanov \cite{Che67}, Rice \cite{Ric68},
Ohtsuka \cite{Oht81}, \cite{Oht85}, \cite{Oht86},
Ohtsuka-Khludnev \cite{O-K00}, and 
Francfort-Marigo \cite{F-M98}.
For more complete list of crack problems and fracture mechanics, 
please see the references in the above papers.

In the Griffith's theory and its various extended theories
such as \cite{F-M98}, the concept of 
the energy release rate $G$ plays an important role.
According to such theories, we treat crack evolutions
in brittle materials with linear elasticity 
under a quasi-static situation,
in which applied boundary loading is supposed to change slowly
and any inertial effect can be ignored.
The elastic energy at a fixed moment is supposed to be  
given by minimization of an elastostatic energy.
According to the Griffith's theory, the surface energy
required in the crack evolution
is supplied by relaxation of the potential energy along crack growth.

Roughly speaking,
the energy release rate $G$ is defined as follows.
Please see the above references for more precise definition.
Let $\Omega_*$ be a bounded domain
in $\R^n$ $(n\geq 2)$, which corresponds to the uncracked
material under consideration.
We assume that a crack $\Sigma$
exist in $\Omega_*$, where $\Sigma$ is the closure of
an $n-1$ dimensional hypersurface.
The cracked elastic body is represented 
by $\Omega_*\setminus\Sigma$.
We consider a virtual crack extension
$\Sigma(t)$ with parameter $t\in [0,T)$,
where 
$\Sigma=\Sigma(0)\subset \Sigma(t_1)\subset\Sigma(t_2)$
$(0\leq ^\forall t_1\leq ^\forall t_2<T)$.

Under the quasi-static assumption,
the elastic potential energy 
$E(t)$ in $\Omega(t):=\Omega_*\setminus\Sigma(t)$
is given by
\begin{equation}\label{energymin}
E(t):=\min_{v}\int_{\Omega_*\setminus\Sigma(t)}
W(x,v(x),\nabla v(x))dx,
\end{equation}
where by $W(x,v(x),\nabla v(x))$ 
the potential energy density including a body force
is denoted, and 
$\min_{v}$ is taken over all possible displacement fields in
$\Omega_*\setminus\Sigma(t)$ with a given boundary
condition. 
For the admissible displacement fields,
a given displacement
field is imposed only on the part $\Gamma_D\subset \partial\Omega_*$.
On the other part $\partial\Omega(t)\setminus \Gamma_D$
including both sides of $\Sigma(t)$,
the normal stress free condition
is imposed for the minimizer on $\partial\Omega(t)$
implicitly.

The energy release rate $G$ at $t=0$
along the virtual crack extension
$\{\Sigma(t)\}_{0\leq t<T}$
is given by
\begin{equation}\label{Glim}
G:=\lim_{t\rightarrow +0}\frac{E(0)-E(t)}
{|\Sigma(t)\setminus\Sigma |}.
\end{equation}
Since $E(t)\leq E(0)$,
$G\geq 0$ follows if the limit exists.
The Griffith's criterion for the brittle crack extension 
is given by $G\geq G_c$, where $G_c$ is an energy required to create
new crack per unit length and it is a constant depending on
the material property and the position.

Cherepanov \cite{Che67} and Rice \cite{Ric68}
studied so-called J-integral for straight crack in two dimensional 
linear elasticity, which is a path-independent integral expressions of the 
energy release rate.
Since these works, theoretical and practical studies of
crack evolutions have been much developed by means of 
such useful mathematical expression of $G$
in two dimensional case.

As an alternative approach to such energy based arguments,
Irwin \cite{Irw47} proposed the notions of
fracture toughness and stress intensity factors
and he developed arguments based on the singularity of 
stress fields.

While most of these mathematically rigorous results
have been restricted to two dimensional linear
elasticity (and often only for straight cracks),
Ohtsuka \cite{Oht81}, \cite{Oht85}, \cite{Oht86}
and Ohtsuka-Khludnev \cite{O-K00} 
developed a mathematical formulation of the
energy release rate for general curved cracks in multi-dimensional
linear or semi-linear elliptic systems.
They proved existence of the energy release rate,
and obtained its expression by a domain integral
and by a generalized J-integral.

Based on the idea in \cite{Oht85}, we shall give a
new mathematical framework
for shape derivative of potential energy
including the energy release rate.
Adopting domain perturbation $\varphi$ of Lipschitz class,
we treat the shape derivative
as an abstract parameter variation problem
in Banach spaces, where $\varphi$ is considered as 
a parameter belonging to a Lipschitz class.
Instead of estimating the limit (\ref{Glim}) directly as in
\cite{Oht85} and \cite{O-K00}, we treat it by means of
the Fr\'{e}chet derivative.

In our approach, the shape derivative of minimum potential energy
is derived as a Fr\'{e}chet derivative in a Banach space
within an abstract parameter variation formulas
and it is given as a domain integral.  
The key tools in the abstract parameter variation setting 
are the implicit function theorem and the Lax-Milgram theorem.

The organization of this paper is as follows.
Abstract parameter variation formulas are
established based on the implicit function theorem
in Banach spaces in Section~\ref{pvfI}.
In Section~\ref{dd}, a framework of Lipschitz deformation of
domains, which includes crack extensions,
is introduced.
Minimization problems with a general potential energy 
in deformed domains are studied in Section~\ref{pemp}
as an application of the abstract parameter variation
formulas in Section~\ref{pvfI}.
Quadratic energy functionals corresponding to
second order linear elliptic equations are treated 
in Section~\ref{alep}. Under a weak regularity assumption,
we show the existence of the shape derivative of the minimum
potential energy and derive its domain integral expression
and a boundary integral formula (J-integral).
The results obtained there
include the results in \cite{Oht85} and \cite{O-K00} under a weaker
assumption for regularity of domain perturbation.
In \cite{O-K00}, they assumed that the domain perturbation
$t\to \varphi(t)$ belongs to $C^2([0,T],W^{2,\infty}(\R^n)^n)$
and derived the domain integral expression,
whereas, in Theorem~\ref{main1}, we prove it under a weaker
assumption $\varphi\in C^1([0,T],W^{1,\infty}(\R^n)^n)$.
For simplicity, description will be made on scalar equation in
Sections \ref{pemp} and \ref{alep}. But our results are
easily extended to elliptic systems such as linear elasticity
problems as shown in \cite{Oht85}.

\subsection{Parameter variation formulas}\label{pvfI} 
\setcounter{equation}{0}

We consider a variational problem with a
parameter in an abstract setting.
For a real valued
functional $J$ defined on  a metric space $S$,
$u_0\in S$ is called a {\it global minimizer} of $J$ in $S$, if
$J(u_0)\leq J(u)$ for all $u\in S$.
If there exists an open set $\cO\subset S$
and $u_0$ is a global minimizer of $J$ in $\cO$,
$u_0\in S$ is called a {\it local minimizer} of $J$.

Let $X$ and $M$ be real Banach spaces.
For open subsets $\U_0\subset X$ and $\cO_0\subset M$,
we consider $J\in C^1(\U_0\times \cO_0,\R)$ and $u\in C^1(\cO_0,\U_0)$.
We assume that $u(\mu)$ is a local minimizer of $J(\cdot,\mu)$ in $\U_0$
for each $\mu\in\cO_0$, and define
$J_*(\mu):=J(u(\mu),\mu)$ for $\mu\in \cO_0$.
Then we have $J_*\in C^1(\cO_0)$ and 
\begin{equation}\label{J*'}
J_*'(\mu)=
D_\mu [J(u(\mu),\mu)]
=\partial_XJ(u(\mu),\mu)[u'(\mu)] 
+
\partial_MJ(u(\mu),\mu)
=
\partial_MJ(u(\mu),\mu),
\end{equation}
where $J_*'$ denotes the Fr\'{e}chet derivative of $J_*$
and $D_\mu$ denotes the Fr\'{e}chet differential operator
with respect to $\mu\in M$.
The symbols $\partial_X$ and $\partial_M$ denote 
the partial Fr\'{e}chet derivative operators for 
$J(u,\mu)$
with respect to $u\in X$ and $\mu\in M$, respectively.
The last equality of (\ref{J*'}) follows from 
$\partial_XJ(u(\mu),\mu)=0\in X'$,
where $X'$ denotes the dual space of $X$.
The formula
\begin{equation}\label{diffmln}
J_*'(\mu)=
\partial_MJ(u(\mu),\mu)
~~~~~
(\mu\in \cO_0),
\end{equation}
is a simple but essential equation in this paper.

The following fundamental theorem states that the formula (\ref{diffmln})
is derived under a weaker assumption for regularity. 
\begin{Th}
\label{naruto}
Let $X$ and $M$ be real Banach spaces.
For $\U_0\subset X$ and an open subset $\cO_0\subset M$,
we consider a real valued functional
$J:\U_0\times \cO_0\rightarrow \R$ and a map 
$u:\cO_0\rightarrow \U_0$. We define
$J_*(\mu):=J(u(\mu),\mu)$ for $\mu\in \cO_0$.
We suppose the following conditions.
{\rm 
\begin{enumerate}
\item {\it $J\in C^0(\U_0\times \cO_0)$,
$J(w,\cdot)\in C^1(\cO_0)$ for $w\in \U_0$, and 
$\partial_MJ\in C^0(\U_0\times\cO_0,M')$.}
\item {\it $u\in C^0(\cO_0,X)$ and 
$u(\mu)$ is a global minimizer of $J(\cdot,\mu)$ in $\U_0$
for each $\mu\in\cO_0$.}
\end{enumerate}
}
\noindent
Then $J_*\in C^1(\cO_0)$ and (\ref{diffmln}) holds. 
\end{Th}
\noindent
{\it Proof.}
We fix $\mu_0\in \cO_0$ and we define $u_0:=u(\mu_0)$ and
\[
r(\mu):=J_*(\mu)-J_*(\mu_0)-\partial_MJ(u_0,\mu_0)[\mu-\mu_0]
~~~~~
(\mu\in\cO_0).
\]
Since $u(\mu)$ is a global minimizer and $u\in C^0(\cO_0,X)$,
if $\mu$ is close to $\mu_0$, we have
\[
r(\mu)\leq J(u_0,\mu)-J(u_0,\mu_0)-\partial_MJ(u_0,\mu_0)[\mu-\mu_0]
=o(\|\mu-\mu_0\|_M),
\]
\begin{eqnarray*}
r(\mu)&\geq&
J(u(\mu),\mu)-J(u(\mu),\mu_0)-\partial_MJ(u_0,\mu_0)[\mu-\mu_0]\\
&=&\int_0^1
\{\partial_MJ(u(\mu),\mu_0+s(\mu-\mu_0))
-\partial_MJ(u_0,\mu_0)\}[\mu-\mu_0]ds
=o(\|\mu-\mu_0\|_M).
\end{eqnarray*}
It follows that $r(\mu)=o(\|\mu-\mu_0\|_M)$ as $\mu\to\mu_0$,
and we obtain the formula (\ref{diffmln}) and $J_*'\in C^0(\cO_0,M')$.
\qed 

\begin{Cor}
Under the condition of Theorem~\ref{naruto}, 
we assume that $\U_0$ is open.
If
$\partial_MJ\in C^k(\U_0\times\cO_0,M')$
and $u\in C^k(\cO_0,X)$, then $J_*\in C^{k+1}(\cO_0)$.
\end{Cor}
\noindent
{\it Proof.}
This immediately follows from the formula (\ref{diffmln}).
\qed

We apply the implicit function theorem in Banach spaces below.
The proof is found in \cite{Ber77} and \cite{Hen81} etc.
For two Banach spaces $X$ and $Y$, $B(X,Y)$ denotes 
the Banach space
which consists of all bounded linear operators
from $X$ to $Y$.
\begin{Th}[Implicit function theorem]
\label{ift}
Let $X$, $Y$, $Z$ be real Banach spaces and 
$U$, $V$ be open sets in $X$ and $Y$, respectively.
We suppose that  
$F:U\times V\rightarrow Z$ and $(x_0,y_0)\in U\times V$
satisfy the conditions;
{\rm
\begin{enumerate}
 \item $F(x_0,y_0)=0$.
 \item $F\in C^0(U\times V,Z)$.
 \item {\it $F(x,\cdot)\in C^1(V,Z)$ for $x\in U$
 and 
 $\partial_YF$ is continuous at $(x,y)=(x_0,y_0)$.}
 \item $(\partial_YF(x_0,y_0))^{-1}\in B(Z,Y)$.
\end{enumerate}}
\noindent
Then there exist a convex open neighborhood of $(x_0,y_0)$,
$U_0\times V_0\subset U\times V$ and $f\in C^0(U_0,V_0)$,
such that, for $(x,y)\in U_0\times V_0$,
$F(x,y)=0$ if and only if $y=f(x)$.
Moreover, if $F\in C^k(U\times V,Z)$ $(k\in \N)$,
then $f\in C^k(U_0,V_0)$.
\end{Th}

From Theorem~\ref{naruto} and the implicit function theorem,
we get the following theorems.
\begin{Th}
\label{fpvfnr}
Let $X$ and $M$ be real Banach spaces
and $\U$ and $\cO$ be open subsets of $X$
and $M$, respectively.
We consider a real valued functional
$J:\U\times \cO\rightarrow \R$ and 
fix $\mu_0\in \cO$. 
We assume
{\rm
\begin{enumerate}
 \item {\it $J(\cdot,\mu)\in C^2(\U)$ 
 for $\mu\in \cO$
and $\partial_XJ\in C^0(\U\times\cO,X')$.}
\item {\it $u_0\in \U$
satisfies $\partial_XJ(u_0,\mu_0)=0$.}
 \item {\it $\partial_X^2J$ is continuous at $(w,\mu)=(u_0,\mu_0)$.}
 \item {\it There exists $\alpha>0$ such that 
 $\partial^2_XJ(u_0,\mu_0)[w,w]\geq \alpha\|w\|_X^2$
 for $w\in X$.}
\end{enumerate}
}\noindent
Then there exist a convex open neighborhood of $(u_0,\mu_0)$,
$\U_0\times\cO_0\subset \U\times\cO$ and 
$u\in C^0(\cO_0,\U_0)$,
such that, for $\mu\in \cO_0$,
the following three conditions are equivalent.\vspace{1ex}\\
~~~{\rm a.} $w\in \U_0$ is a local minimizer of $J(\cdot,\mu)$\vspace{1ex}\\
~~~{\rm b.} $w\in \U_0$ satisfies $\partial_XJ(w,\mu)=0$.\vspace{1ex}\\
~~~{\rm c.} $w=u(\mu)$.\\
In this case, $u(\mu)$ is a global minimizer of $J(\cdot,\mu)$ on $\U_0$.
\end{Th}
\noindent
{\it Proof.}
We define a map $F:=\partial_XJ$ 
from $\U\times\cO$ to $X'$
and apply Theorem~\ref{ift}
at $(w,\mu)=(u_0,\mu_0)$.
From assumption 4 and the Lax-Milgram theorem,
$\partial_XF(u_0,\mu_0)=\partial_X^2J(u_0,\mu_0)$ becomes
a linear topological
isomorphism from $X$ to $X'$.
Then, from the implicit function theorem 
there exist a convex open neighborhood of $(u_0,\mu_0)$,
$\U_0\times\cO_0\subset \U\times\cO$ and 
$u\in C^0(\cO_0,U_0)$,
such that, for $\mu\in\cO_0$, $w\in \U_0$ satisfies
$\partial_XJ(w,\mu)=0$
if and only if $w=u(\mu)$.

From the continuity of $\partial^2_XJ$
at $(u_0,\mu_0)$, without loss of generality,  
(after replacing $\U_0$ and $\cO_0$ with smaller ones if needed)
we can assume that
\begin{equation}\label{alpha2}
\partial^2_XJ(v,\mu)[w,w]\geq \frac{\alpha}{2}\|w\|_X^2
~~~~~(^\forall w\in X,~^\forall (v,\mu)\in \U_0\times\cO_0).
\end{equation}

For $\mu\in \cO_0$, if $w\in \U_0$ is a local
minimizer of $J(\cdot,\mu)$ in $\U_0$,
the $\partial_XJ(w,\mu)=0$ follows.
Conversely, if $w\in \U_0$ satisfies $\partial_XJ(w,\mu)=0$,
$w$ is a local minimizer in $\U_0$ from the condition (\ref{alpha2}).
It also follows from  (\ref{alpha2}) that
$u(\mu)$ is a global minimizer of $J(\cdot,\mu)$ in $\U_0$.
\qed

\begin{Th}\label{addth}
Under the condition of Theorem~\ref{fpvfnr},
we additionally assume that 
$\partial_XJ\in C^k(\U\times\cO,X')$
for some $k\in \N$.
Then $u\in C^k(\cO_0,\U_0)$ holds.
\end{Th}
\noindent
{\it Proof.}
The assertion follows from the implicit function theorem.
\qed 

Under the condition of Theorem~\ref{fpvfnr},
we define 
\[
J_*(\mu):=J(u(\mu),\mu)
~~~~~(\mu\in\cO_0).
\]
As a consequence of Theorem~\ref{addth}, a sufficient condition 
for $J_*\in C^1(\cO_0)$ is $J\in C^1(\U\times\cO)$ and 
$\partial_XJ\in C^1(\U\times\cO,X')$. However, the condition
$\partial_XJ\in C^1(\U\times\cO,X')$ is not necessary due to
Theorem~\ref{naruto}.

\begin{Th}
\label{ckcondition}
Under the condition of Theorem~\ref{fpvfnr},
we additionally assume that
$J\in C^k(\U\times\cO)$ for some $k\in\N$,
then $J_*\in C^k(\cO_0)$
and it satisfies {\rm (\ref{diffmln})}.
\end{Th}
\noindent
{\it Proof.}
From Theorem~\ref{naruto}, $J_*\in C^1(\cO_0)$
and (\ref{diffmln}) immediately follows.
Since $u\in C^{k-1}(\cO_0,X)$ follows from Theorem~\ref{addth},
$J_*\in C^k(\cO_0)$ is obtained from the formula (\ref{diffmln}).
\qed

Let us consider the case $k=1$ in Theorem~\ref{ckcondition},
where $J_*\in C^1(\cO_0)$ is derived under the conditions
$J\in C^1(\U\times\cO)$ and $J(\cdot,\mu)\in C^2(\U)$.
In this case, $u\in C^0(\cO_0,\U_0)$ holds from
Theorem~\ref{fpvfnr} but $u\not\in C^1(\cO_0,\U_0)$
in general.
In order to obtain $u\in C^1(\cO_0,\U_0)$, we need to
assume $\partial_XJ\in C^1(\U\times\cO,X')$
(Theorem~\ref{addth}).
We have H\"{o}lder regularity of $u$ under the condition
of Theorem~\ref{ckcondition} with $k=1$.
\begin{Prop}
\label{uholder}
Under the condition of Theorem~\ref{fpvfnr},
we additionally assume that
$J\in C^1(\U\times\cO)$,
then we have
\[
\|u(\mu)-u_0\|_X=o\left(\|\mu-\mu_0\|_M^{1/2}\right)
~~~~~
\mbox{as}~~\|\mu-\mu_0\|_M\rightarrow 0.
\]
\end{Prop}
\noindent
{\it Proof.}
From the proof of Theorem~\ref{ift}, 
there exists $C>0$ such that
\begin{equation}\label{fromproof}
\|u(\mu)-u_0\|_X\leq
C\|\partial_XJ(u_0,\mu)\|_{X'}
~~~
(\mu\in \cO_0).
\end{equation}
Let $\rho_0>0$ with 
$\{v\in X;~\|v-u_0\|_X\leq \rho_0\}\subset\U_0$.
For $h\in X$ with $\|h\|_X=1$, $\mu\in\cO_0$ and $\rho\in (0,\rho_0]$,
we have
\[
J(u_0+\rho h,\mu)=J(u_0,\mu)+\rho\partial_XJ(u_0,\mu)[h]
+\rho^2\int_0^1(1-s)\partial_X^2J(u_0+s\rho h,\mu)[h,h]ds.
\]
\begin{eqnarray*}
\lefteqn{\partial_XJ(u_0,\mu)[h]} \\
&=&
\partial_XJ(u_0,\mu)[h]-\partial_XJ(u_0,\mu_0)[h]\\
&=&
\left\{
\rho^{-1}(J(u_0+\rho h,\mu)-J(u_0,\mu))
-\rho\int_0^1(1-s)\partial_X^2J(u_0+s\rho h,\mu)[h,h]ds
\right\}\\
&&
-\left\{
\rho^{-1}(J(u_0+\rho h,\mu_0)-J(u_0,\mu_0))
-\rho\int_0^1(1-s)\partial_X^2J(u_0+s\rho h,\mu_0)[h,h]ds
\right\}\\
&=&
\rho^{-1}\int_0^1
\{\partial_MJ(u_0+\rho h,\mu_0+t(\mu-\mu_0))
-\partial_MJ(u_0,\mu_0+t(\mu-\mu_0))\}[\mu-\mu_0]dt\\
&&
-\rho\int_0^1(1-s)
\{\partial_X^2J(u_0+s\rho h,\mu)-\partial_X^2J(u_0+s\rho h,\mu_0)\}[h,h]ds
\end{eqnarray*}
For $r>0$, we define
\[
S(r):=\left\{(w,\lambda)\in X\times M;~
\|w-u_0\|_X\leq r,\|\lambda-\mu_0\|_M\leq r^2\right\},
\]
\begin{eqnarray*}
\omega(r)&:=&
\sup_{(w,\lambda)\in S(r)}
\|\partial_MJ(w,\lambda)-\partial_MJ(u_0,\mu_0)\|_{M'}\\
&&+
\sup_{(w,\lambda)\in S(r)}
\|\partial_X^2J(w,\lambda)-\partial_X^2J(u_0,\mu_0)\|_{B_2(X,\R)}.
\end{eqnarray*}
We remark that $\omega(r)\rightarrow 0$ as $r\rightarrow +0$.
Choosing
$\rho:=\|\mu-\mu_0\|_M^{1/2}$, we obtain
\[
\|\partial_XJ(u_0,\mu)\|_{X'}\leq 2\omega(\rho)\rho
~~~~~~~
(\mu\in \cO_0,~\|\mu-\mu_0\|_M\leq \rho_0^2).
\]
Hence, from (\ref{fromproof}), we have
\[
\|u(\mu)-u_0\|_X\leq
C\|\partial_XJ(u_0,\mu)\|_{X'}
\leq
2C\omega(\rho)\rho=o(\rho).
\]
\qed

Under the conditions of Theorem~\ref{fpvfnr},
$\partial_X^2J(u(\mu),\mu)$ can be regarded as
a linear topological isomorphism from $X$ to $X'$ 
from the Lax-Milgram theorem.
Therefore, we can define $\Lambda (\mu)\in B(X',X)$
which satisfies
\[
\partial_X^2J(u(\mu),\mu)[\Lambda (\mu)h,w]
=h[w]~~~~~(^\forall w\in X,~^\forall h\in X').
\]
The Fr\'{e}chet derivative of the local minimizer $u(\mu)$
with respect to parameter $\mu$ is given by the next proposition.
\begin{Prop}\label{h0prop}
Under the condition of Theorem~\ref{addth} with $k=1$,
\begin{equation}\label{dotu}
u'(\mu)=-\Lambda(\mu)h_0(\mu)~~~~~
(\mu\in\cO_0),
\end{equation}
holds, where $h_0(\mu):=\partial_M\partial_X J(u(\mu),\mu)\in B(M,X')$.
\end{Prop}
\noindent
{\it Proof.}
Differentiating $\partial_XJ(u(\mu),\mu)=0\in X'$ by $\mu$, we have
\[
\partial_X^2J(u(\mu),\mu)[u'(\mu)]
+
\partial_M\partial_XJ(u(\mu),\mu)
=0\in B(M,X').
\]
This is equivalent to 
(\ref{dotu}) from the Lax-Milgram theorem.
\qed 

\begin{Prop}\label{J*2}
Under the condition of Theorem~\ref{fpvfnr},
we additionally assume that
$J\in C^2(\U\times\cO)$
then $J_*\in C^2(\cO_0)$
and it satisfies
\[
J_*''(\mu)[\mu_1,\mu_2]
=\partial_M^2J(u(\mu),\mu)[\mu_1,\mu_2]
-~_{X}\hspace{-0.1ex}\langle \Lambda(\mu)h_0(\mu)[\mu_1],
~h_0(\mu)[\mu_2] \rangle_{X'},
\]
for
$\mu\in\cO_0$ and $\mu_1,\mu_2\in M$.
\end{Prop}
\noindent
{\it Proof.}
Differentiating the formula (\ref{diffmln}) by $\mu$
and substituting (\ref{dotu}),
we obtain the formula.
\qed
\subsection{Lipschitz deformation of domains}\label{dd} 
\setcounter{equation}{0}
We study a domain deformation with Lipschitz transform 
$\varphi:\Omega\rightarrow \varphi(\Omega)$, where $\Omega$ is a
bounded domain in $\R^n$ $(n\in\N)$ 
and $\varphi$ is a $\R^n$-valued Lipschitz
function.
The identity map on $\R^n$ is denoted by
$\varphi_0(x)=x$ ($x\in\R^n$). 

For a function $u:\Omega\rightarrow \R^k$, 
we define
\[
|u|_{\mbox{\scriptsize Lip},\Omega}:=\sup_{x,y\in\Omega,x\neq y}
\frac{|u(x)-u(y)|}{|x-y|},
\]
where $|\cdot |$ denotes the Euclidean norm in $\R^n$ or $\R^k$.
If $|u|_{\mbox{\scriptsize Lip},\Omega}<\infty$, $u$ is called 
uniformly Lipschitz 
continuous
on $\Omega$.
It is known that, for $u\in W^{1,\infty}(\Omega)$, 
there is $\tilde{u}\in C^0(\Omega)$
such that $\tilde{u}(x)=u(x)$ a.e. $x\in\Omega$,
in other words, we can regard $W^{1,\infty}(\Omega)\subset C^0(\Omega)$.
If $\Omega$ is convex,
$W^{1,\infty}(\Omega)=C^{0,1}(\ov{\Omega})$
as a subset of $C^0(\Omega)$. Moreover, if $k=1$, we have
\[
\|\nabla u\|_{L^\infty(\Omega)}=
|u|_{\mbox{\rm\scriptsize Lip},\Omega}
~~~~~
(u\in W^{1,\infty}(\Omega)\cap C^0(\Omega)).
\]
In the following argument, we fix a bounded convex domain
$\Omega_0\subset\R^n$ ($n\geq 2$),
and we identify 
$W^{1,\infty}(\Omega_0,\R^n)$
with $C^{0,1}(\ov{\Omega_0},\R^n)$.

\begin{Prop}\label{bilipschitz}
Suppose that $\varphi\in W^{1,\infty}(\Omega_0,\R^n)$
satisfies 
$|\varphi-\varphi_0|_{\mbox{\rm\scriptsize Lip},\Omega_0}<1$. 
Then $\varphi$ is a
bi-Lipschitz transform from $\Omega_0$ to $\varphi(\Omega_0)$,
i.e. 
$\varphi$ is bijective from $\Omega_0$ onto an open set
and $\varphi$ and $\varphi^{-1}$ are both uniformly Lipschitz continuous.
Moreover, we have
\begin{equation}\label{det-positive}
\mbox{\rm ess-}\inf_{\Omega_0}\left(\det\nabla \varphi^\srT \right)
\geq \left( 1- |\varphi-\varphi_0|_{\mbox{\rm\scriptsize Lip},\Omega_0}
\right)^n >0.
\end{equation}
where $\nabla \varphi^\srT$ is the Jacobian matrix defined by
\[
\nabla \varphi^\srT (x):=
\left(\frac{\partial\varphi_j}{\partial x_i}(x)\right)
_{i\downarrow ,j\rightarrow }\in\R^{n\times n}
~~~~~
(x=(x_1,\cdots,x_n)^\srT \in \Omega_0).
\]
\end{Prop}
\noindent
{\it Proof.}
Let $\mu:=\varphi-\varphi_0$
and $\theta:=|\mu|_{\mbox{\rm\scriptsize Lip},\Omega_0}\in (0,1)$.
First, we show that $\varphi(\Omega_0)$ is open.
We arbitrarily fix $x_0\in\Omega_0$
and define $y_0:=\varphi(x_0)$. Let 
$\delta >0$ 
such that $\ov{B_\delta(x_0)}\subset \Omega_0$,
where $B_\delta(x_0):=\{x\in\R^n;~|x-x_0|<\delta\}$.
For $y\in B_{(1-\theta)\delta}(y_0)$, we show that 
$y\in \varphi(\Omega_0)$.
It is easily checked that $T(\xi):=y-\mu(\xi)$ is a uniform
contraction on $\ov{B_\delta(x_0)}$.
From the contraction mapping theorem,
there is a fixed point $x=T(x)=y-\mu(x)$ in $\ov{B_\delta(x_0)}$,
that is $y=\varphi(x)$.
Hence, $\varphi(\Omega_0)$ is an open set.
Since
$|\varphi(x_1)-\varphi(x_2)|\geq |x_1-x_2|-|\mu(x_1)-\mu(x_2)|
\geq (1-\theta)|x_1-x_2|$
for $x_1,x_2\in \Omega_0$,
it follows that $\varphi$ is injective and $\varphi^{-1}$
satisfies uniform Lipschitz condition on $\varphi(\Omega_0)$.

From Rademacher's theorem (see \cite{E-G92}, \cite{Zie89}),
$\mu$ is differentiable almost everywhere and the derivative 
coincides with the distributional derivative almost everywhere, i.e.,
there exists $\cN\subset \Omega_0$ with $\cL^n(\cN)=0$ such that
\[
\nabla^\srT \mu(x)y=
\lim_{h\to 0}
\frac{\mu(x+hy)-\mu(x)}{h}
~~~~~
(x\in\Omega_0\setminus\cN),
\]
where $\nabla^\srT \mu=(\nabla \mu^\srT)^\srT$.
It follows that
\[
\left| \nabla^\srT \mu(x)y \right| 
\leq \theta \,|y|
~~~~~
(x\in\Omega_0\setminus\cN),
\]
and that the moduli of all the eigenvalues of $\nabla^\srT \mu(x)$ 
for $x\in\Omega_0\setminus\cN$
are bounded by $\theta$. Hence, we obtain
\[
\det\left(\nabla^\srT \varphi (x)\right)
=
\det\left(I+\nabla^\srT \mu (x)\right)
\geq (1-\theta)^n
~~~~~
(x\in\Omega_0\setminus\cN).
\]
\qed

We fix an open set $\Omega$ which satisfies 
$\ov{\Omega}\subset\Omega_0$ and
the deformed domain $\varphi(\Omega)$
is denoted by $\Omega(\varphi)$
under the condition of Proposition~\ref{bilipschitz}, hereafter.
We define a pushforward operator $\varphi_*$
which transforms
a function $v$ on $\Omega$ to a function 
$\varphi_*v:=v\circ\varphi^{-1}$ 
on $\Omega(\varphi)$,
if $\varphi$ satisfies Proposition~\ref{bilipschitz}.
We define 

\[
A(\varphi):=\left(\nabla \varphi^\srT \right)^{-1}
\in L^\infty(\Omega_0,\R^{n\times n}),
~~~
\kappa(\varphi):=\det\nabla \varphi^\srT 
\in L^\infty(\Omega_0,\R).
\]
These Jacobi matrices and Jacobian appear in the pullback of
differentiation and integration on $\Omega(\varphi)$ to $\Omega$.
For a function $v$ on $\Omega$, 
we have
\begin{equation}\label{form1}
[\nabla(\varphi_*v)]\circ \varphi
=A(\varphi)\nabla v
~~~
\mbox{a.e in}~\Omega
~~~
(v\in W^{1,1}(\Omega)),
\end{equation}
\begin{equation}\label{form2}
\int_{\Omega(\varphi)}
(\varphi_*v)(y)dy
=
\int_\Omega v(x) \kappa(\varphi)(x)dx
~~~~~
(v\in L^1(\Omega))
\end{equation}
These equalities are well known in the case $\varphi\in C^1$.
However, for $\varphi\in C^{0,1}$, these are not so trivial.
See, \cite{E-G92} and \cite{Zie89} etc. for details.
We omit the proof of the next proposition
since it is clear from (\ref{form1}) and (\ref{form2}).
\begin{Prop}\label{mM}
Under the condition of Proposition~\ref{bilipschitz},
for $p\in [1,\infty]$, $\varphi_*$ is a linear
topological isomorphism from $L^p(\Omega)$
onto $L^p(\Omega(\varphi))$, and a linear
topological isomorphism from $W^{1,p}(\Omega)$
onto $W^{1,p}(\Omega(\varphi))$.
\end{Prop}
The following theorem plays an essential role in
the application to the shape derivatives.
\begin{Th}\label{kappaA}
Let $\Omega$ be an open subset of $\Omega_0$.
{\rm
\begin{enumerate}
\item
{\it $\kappa\in C^\infty(W^{1,\infty}(\Omega,\R^n),L^\infty(\Omega))$,
and $\kappa'(\varphi_0)[\mu]=\mbox{\rm div}\mu$
for $\mu\in W^{1,\infty}(\Omega,\R^n)$.}
\item
{\it 
$A\in C^\infty(\cO_0(\Omega),L^\infty(\Omega,\R^{n\times n}))$,
where 
\[
\cO_0(\Omega):=\{\varphi\in W^{1,\infty}(\Omega,\R^n);~
\mbox{\rm ess-}\inf_{\Omega} \kappa(\varphi)>0\}.
\]
In particular, $A'(\varphi_0)[\mu]=-\nabla\mu^\srT $ holds
for $\mu\in W^{1,\infty}(\Omega,\R^n)$.}
\end{enumerate}
}
\end{Th}
\noindent
{\it Proof.}
Since the determinant is a polynomial of degree $n$,
it is clear that 
$\kappa$ belongs to $C^\infty(W^{1,\infty}(\Omega,\R^n),L^\infty(\Omega))$.
For fixed $\mu\in W^{1,\infty}(\Omega,\R^n)$, we define
\[
m_{ij}(t):=\delta_{ij}+t\frac{\partial \mu_j}{\partial x_i}
\in L^\infty(\Omega)
~~~~~
(i,j=1,\cdots,n,~t\in\R),
\]
where $\delta_{ij}$ is the Kronecker's delta.
Then we have
\begin{eqnarray*}
\kappa'(\varphi_0)[\mu]
&=&
\left.\frac{d}{dt}\right|_{t=0}
\kappa(\varphi_0+t\mu)
=
\left.\frac{d}{dt}\right|_{t=0}
\det\left( m_{ij}(t) \right)\\
&=&
\left.\frac{d}{dt}\right|_{t=0}
\left(
\sum_{\sigma\in S_n}
\mbox{\rm sgn}(\sigma)
m_{1\sigma(1)}(t)\cdots m_{n\sigma(n)}(t) 
\right)\\
&=&
\left.\frac{d}{dt}\right|_{t=0}
\left(
m_{11}(t)\cdots m_{nn}(t) 
\right)
=
\sum_{i=1}^n
m_{11}(0)\cdots m_{ii}'(0)\cdots m_{nn}(0) 
=
\mbox{\rm div}\mu.
\end{eqnarray*}

Let the $(i,j)$ component of $A(\varphi)$ be denoted by 
$a_{ij}(\varphi)\in L^\infty(\Omega)$. Then we have
$a_{ij}(\varphi)=\alpha_{ij}(\varphi)/\kappa(\varphi)$,
where $\alpha_{ij}(\varphi)$ is the $(i,j)$ cofactor of
$\nabla \varphi^\srT $, which is a polynomial of 
$\frac{\partial \varphi_k}{\partial x_l}$ of degree $n-1$.
Since $\mbox{\rm ess-}\inf_{\Omega} \kappa(\varphi)>0$
for $\varphi\in\cO_0(\Omega)$, 
$a_{ij}\in C^\infty(\cO_0(\Omega),L^\infty(\Omega))$
follows.
For fixed $\mu\in W^{1,\infty}(\Omega,\R^n)$,
differentiating the identity
\[
A(\varphi_0+t\mu)(I+t\nabla\mu^\srT )=I
~~~~~
\mbox{
($I$: identity matrix of degree $n$)},
\]
by $t\in\R$ at $t=0$, we have
\[
A'(\varphi_0)[\mu]+A(\varphi_0)\nabla\mu^\srT =O.
\]
Since $A(\varphi_0)=I$, we have 
$A'(\varphi_0)[\mu]=-\nabla\mu^\srT $.
\qed 

We define an open subset of $W^{1,\infty}(\Omega_0,\R^n)$
as
\begin{equation}\label{OOmega}
\cO(\Omega):=\left\{\varphi\in W^{1,\infty}(\Omega_0,\R^n);~
|\varphi-\varphi_0|_{\mbox{\rm\scriptsize Lip},\Omega_0}<1,~
\ov{\Omega(\varphi)}\subset\Omega_0
\right\}.
\end{equation}
\begin{Prop}\label{phimu}
We assume $\mu\in W^{1,\infty}(\Omega_0,\R^n)$ 
with $\supp (\mu)\subset \Omega$.
{\rm
\begin{enumerate}
 \item {\it If $|\mu|_{\mbox{\rm\scriptsize Lip},\Omega_0}<1$,
 then $\varphi=\varphi_0+\mu$ is a bi-Lipschitz transform
 from $\Omega$ onto itself.}
 \item {\it For $t\in\R$ with $|t\mu|_{\mbox{\rm\scriptsize Lip},\Omega_0}<1$,
 we define a bi-Lipschitz transform 
 $\varphi(t)=\varphi_0+t\mu$ 
 from $\Omega$ to itself. 
 Let $l\in \{0,1\}$ and $p\in [1,\infty]$.
 Suppose that $f\in W^{l,p}(\Omega)$ if
$p\in [1,\infty)$, and 
$f\in C^{l}(\Omega)\cap W^{l,\infty}(\Omega)$ if
$p=\infty$. Then $\varphi(t)_*f\rightarrow f$
strongly in $W^{l,p}(\Omega)$ as $t\rightarrow 0$. }
\end{enumerate}
}
\end{Prop}
\noindent
{\it Proof.}
From Proposition~\ref{bilipschitz}, 
claim 1 is clear. For claim 2, 
let us  fix $t_0>0$ with 
$|t_0\mu|_{\mbox{\rm\scriptsize Lip},\Omega_0}<1$.
Then, from Proposition~\ref{mM}, 
there exist $C>0$ such that
the following inequalities hold for $|t|\leq t_0$,
\[
\|\varphi(t)_*f-f\|_{W^{l,p}(\Omega)} 
=
\|\varphi(t)_*(f-f\circ \varphi(t))\|_{W^{l,p}(\Omega)}\leq
C\|f-f\circ \varphi(t)\|_{W^{l,p}(\Omega)}.
\]
Since $[\varphi\mapsto f\circ\varphi]\in 
C^0(\cO(\Omega),W^{l,p}(\Omega))$, we obtain 
\[
\|f-f\circ \varphi(t)\|_{W^{l,p}(\Omega)}=
\|f\circ\varphi_0-f\circ \varphi(t)\|_{W^{l,p}(\Omega)}
\rightarrow 0,
\]
as $t\rightarrow 0$.
\qed
\subsection{Potential energy in deformed domains}\label{pemp} 
\setcounter{equation}{0}
In this section, we consider
minimization problems of an abstract potential energy
in deformed domains.
Some concrete examples in linear elliptic equations
will be given in Section~\ref{alep}.

Let $\Omega_0$ be a fixed bounded convex open set of $\R^n$ ($n\in\N$).
We consider an open set $\Omega$ whose closure is contained in $\Omega_0$.
For $v\in H^1(\Omega)=W^{1,2}(\Omega)$, we introduce
the following energy functional:
\[
E(v,\Omega):=\int_{\Omega}
W(x,v(x),\nabla v(x))dx,
\]
where 
\[
W(\xi,\eta,\zeta)\in \R~~\mbox{for}~~
(\xi,\eta,\zeta)\in\Omega_0\times \R\times \R^n,
\]
is a given energy density function.
We assume some suitable regularity conditions
and boundedness of its derivatives in the following argument.
For simplicity, the partial derivatives of $W$ with respect to
$\xi$, $\eta$ and $\zeta$ will be denoted by
$\nabla_\xi W=(\frac{\partial W}{\partial \xi_1},
\cdots ,\frac{\partial W}{\partial \xi_n})^\srT $,
$W_\eta=\frac{\partial W}{\partial \eta}$, and
$\nabla_\zeta W=(\frac{\partial W}{\partial \zeta_1},
\cdots ,\frac{\partial W}{\partial \zeta_n})^\srT $,
respectively. Moreover, for $v\in H^1(\Omega)$, we often
write $W(v(x))=W(x,v(x),\nabla v(x))$,
$\nabla_\xi W(v(x))=\nabla_\xi W(x,v(x),\nabla v(x))$, etc.

We consider the following minimization problem.
\begin{Prob}\label{probomega}
Let $V$ be a closed subspace of $H^1(\Omega)$
with $H^1_0(\Omega)\subset V\subset H^1(\Omega)$,
and let $V(g):=\{v\in H^1(\Omega);~v-g\in V\}$
for $g\in H^1(\Omega)$.
For given $g\in H^1(\Omega)$, find a local minimizer 
$u$ of $E(\cdot,\Omega)$ in $V(g)$, i.e.
$u\in V(g)$ and there exists $\rho>0$ such that
\begin{equation}\label{minimizing}
E(u,\Omega)\leq E(w,\Omega)
~~~~~
(^\forall w\in V(g)~\mbox{\rm with}~\|w-u\|_{H^1(\Omega)}<\rho).
\end{equation}
\end{Prob}

If $u$ is a local minimizer, under suitable regularity conditions 
for $W$,
formally we obtain the following variation formula:
\begin{equation}\label{uvar}
\int_\Omega
\left\{
W_\eta(u(x))v(x)
+
\nabla_\zeta W(u(x))\cdot\nabla v(x)
\right\}dx=0
~~~~~(^\forall v\in V\subset H^1_0(\Omega)).
\end{equation}
This implies
\[
-\mbox{div}\left[\nabla_\zeta W(u(x))\right]
+W_\eta(u(x))
=0~~~~~
\mbox{in}~\Omega .
\]

For fixed $\Omega$ and $V\subset H^1(\Omega)$ as Problem~\ref{probomega},
we consider a family of minimization problem parametrized by
$\varphi\in \cO(\Omega)$, where $\cO(\Omega)$ is 
defined by (\ref{OOmega}).
We define an affine space in $H^1(\Omega(\varphi))$:
\[
V(\varphi,g):=\varphi_*(V(g))=
\{v\in H^1(\Omega(\varphi));~
\varphi_*^{-1}(v)-g\in V \}
~~~~~
(\varphi\in \cO(\Omega),~g\in H^1(\Omega)).
\]
\begin{Prob}\label{probophi}
For given $\varphi\in \cO(\Omega)$ and 
$g\in H^1(\Omega)$, find a local minimizer 
$u(\varphi)$ of $E(\cdot,\Omega(\varphi))$ in $V(\varphi,g)$, 
i.e.
$u(\varphi)\in V(\varphi,g)$ and there exists $\rho>0$ such that
\begin{equation}\label{mini2}
E(u(\varphi),\Omega(\varphi))\leq E(w,\Omega(\varphi))
~~~~~
(^\forall w\in V(\varphi,g)~\mbox{\rm with}
~\|w-u(\varphi)\|_{H^1(\Omega(\varphi))}<\rho).
\end{equation}
\end{Prob}
We define
\begin{equation}\label{E*}
E_*(\varphi):=
E(u(\varphi),\Omega(\varphi)),
\end{equation}
for a local minimizer $u_*(\varphi)$.
Using the formulas (\ref{form1}) and (\ref{form2}),
we define
\begin{equation}\label{EE}
\E(v,\varphi):=
\int_\Omega
W\left(
\varphi(x),~v(x),~[A(\varphi)(x)]\nabla v(x)
\right)\kappa(\varphi)(x)dx~~~~~
(v\in H^1(\Omega), ~\varphi\in\cO(\Omega)
),
\end{equation}
then we have
\[
E(\varphi_*v,\Omega(\varphi))=\E(v,\varphi)~~~~~
(v\in H^1(\Omega)).
\]
The following theorem is a direct consequence of Theorem~\ref{naruto}.
\begin{Th}\label{direct}
Suppose that $\E$ defined by {\rm (\ref{EE})}
belongs to $C^0(H^1(\Omega)\times\cO(\Omega))$
and that $\E(v,\cdot)\in C^1(\cO(\Omega))$
for $v\in H^1(\Omega)$ and 
$\partial_\varphi\E\in 
C^0(H^1(\Omega)\times\cO(\Omega),~(W^{1,\infty}(\Omega_0,\R^n))')$,
where $\partial_\varphi\E$ denotes the partial Fr\'{e}chet 
derivative of $\E$ 
with respect to $\varphi\in
\cO_0\subset W^{1,\infty}(\Omega_0,\R^n)$.
Let $\U_0\subset V$ and
$\cO_0\in\cO(\Omega)$ be open subsets with $\varphi_0\in \cO_0$,
and we define 
$\U(\varphi,g):=\varphi_*(\U_0+g)\subset V(\varphi,g)$.
If $u(\cdot)\in C^0(\cO_0,H^1(\Omega))$ and 
$u(\varphi)$ is a global minimizer of 
$E(\cdot,\Omega(\varphi))$ in $\U(\varphi,g)$
for all $\varphi\in \cO_0$, 
then we have $E_*\in C^1(\cO_0)$ and
\[
E_*'(\varphi)=\partial_\varphi\E(\varphi_*^{-1}u(\varphi),\varphi).
\]
In particular, we have
\[
E_*'(\varphi_0)=\partial_\varphi\E(u(\varphi_0),\varphi_0).
\]
\end{Th}
\noindent
{\it Proof.}
Let
$v(\varphi):=\varphi_*^{-1}u(\varphi)-g\in V$ and
$\bar{\E}(v,\varphi):=\E(v+g,\varphi)$.
Then $u(\varphi)$ is a local minimizer of $E(\cdot,\Omega(\varphi))$
in $V(\varphi,g)$, if and only if
$v(\varphi)$ is a local minimizer of $\bar{\E}(\cdot,\varphi)$ in $V$.
We apply Theorem~\ref{naruto} to $\bar{\E}$
and obtain $E_*\in C^1(\cO_0)$ and 
\[
E_*'(\varphi)=\partial_\varphi\bar{\E}(v(\varphi),\varphi)
=\partial_\varphi\E(v(\varphi)+g,\varphi)
=\partial_\varphi\E(\varphi_*^{-1}u(\varphi),\varphi).
\]
\qed

The minimizer $u(\varphi_0)$ is denoted by $u$ hereafter.
Under the suitable regularity conditions for $W(\xi,\eta,\zeta)$,
for $\mu\in W^{1,\infty}(\Omega_0,\R^n)$,
we have
\begin{eqnarray}
\lefteqn{\partial_\varphi\E (u,\varphi_0)[\mu]} \\
&=&
\frac{d}{dt}\left.\int_\Omega
W\left(
x+t\mu(x),~u(x),~[A(\varphi_0+t\mu)(x)]
\nabla u(x)
\right)\kappa(\varphi_0+t\mu)(x)dx\right|_{t=0}\nonumber\\
&=&
\int_\Omega
\left(
\nabla_\xi W(u)\cdot \mu 
-(\nabla_\zeta W(u))^\srT (\nabla\mu^\srT )\nabla u
+W(u)\DIV \mu
\right)dx.\label{dEformula}
\end{eqnarray}

Using the above formula, 
we also consider the inner variation, which is
another type of variation for
Problem~\ref{probomega},
as follows.
\begin{Th}
Under the conditions of Theorem~\ref{direct},
if $u=u(\varphi_0)$ is a local minimizer of Problem~\ref{probomega}, 
then we have
\[
\partial_\varphi\E(u,\varphi_0)[\mu]=0
~~~~~
(\mu\in W^{1,\infty}(\Omega),~
\supp (\mu)\subset\Omega).
\]
\end{Th}
\noindent
{\it Proof.}
We suppose that 
$\mu\in W^{1,\infty}(\Omega)$ with $\supp (\mu)\subset\Omega$.
We define $\varphi(t)=\varphi_0+t\mu$ for $t\in\R$.
From Proposition~\ref{mM} and \ref{phimu}, if 
$|t\mu|_{\mbox{\rm\scriptsize Lip},\Omega_0}<1$,
the corresponding pushforward operator $\varphi(t)_*$ is
a linear topological isomorphism from $H^1(\Omega)$ onto itself.
Moreover, from the formulas (\ref{form1}) and 
(\ref{form2}),
\[
\lim_{t\rightarrow 0}\|\varphi(t)_*u-u\|_{H^1(\Omega)}=0,
\]
holds.
Since
\[
E(\varphi(t)_*u,\Omega)\geq E(u,\Omega)=E(\varphi(0)_*u,\Omega)
~~~~~
(|t|<|\mu|_{\mbox{\rm\scriptsize Lip},\Omega_0}^{-1}),
\]
we obtain
\[
0=\left.\frac{d}{dt}E(\varphi(t)_*u,\Omega)\right|_{t=0}
=\left.\frac{d}{dt}\E(u,\varphi(t))\right|_{t=0}
=\partial_\varphi\E(u,\varphi_0)[\mu].
\]
\qed

\subsection{Application to linear elliptic problems}\label{alep} 
\setcounter{equation}{0}

In this section, we consider a second order linear elliptic
equation with a quadratic potential energy
including the anti-plane displacement model of
two dimensional linear elasticity.

\begin{Ex}\label{linex}
{\rm
We consider the following potential energy.
\begin{equation}\label{linearprob1}
W(\xi,\eta,\zeta)=\half
\left(
\zeta^\srT B(\xi)\zeta+b(\xi)\eta^2
\right)-f(\xi)\eta,
\end{equation}
\[
E(v,\Omega)=\int_\Omega
\left\{
\half
\left(
\nabla^\srT\! v\, B(x)\nabla v+b(x)v^2
\right)-f(x)v
\right\}dx,
\]
where 
$k\in\N\cup\{0\}$ and $B(\xi)$ is an $n\times n$ symmetric matrix
(which is denoted by $\R^{n\times n}_{\mbox{\,\scriptsize\rm sym}}$)
and it satisfies 
\[
B\in C^k(\Omega_0,\R^{n\times n}_{\mbox{\,\scriptsize\rm sym}})
~~~\mbox{and}~~~
~^\exists \beta_0>0~\mbox{s.t.}~
\zeta^\srT B(\xi)\zeta\geq \beta_0|\zeta|^2
~~~~~
(^\forall \xi\in \Omega_0,~^\forall \zeta\in \R^n),
\]
and
\begin{equation}\label{bfcond1}
b\in C^k(\Omega_0,\R),
~~~~~
f\in W^{k,2}(\Omega_0,\R),
\end{equation}
with the condition $b(\xi)\geq 0$ for $\xi\in \Omega_0$.
We remark that
\begin{eqnarray*}
\nabla_\xi W(\xi,\eta,\zeta)
&=&
\half
\left(
\sum_{i,j=1}^n \nabla b_{ij}(\xi)\zeta_i \zeta_j 
+\nabla b(\xi)\eta^2
\right)-\nabla f(\xi)\eta, \\~\\
W_\eta (\xi,\eta,\zeta)
&=&
b(\xi)\eta-f(\xi), \\~\\
\nabla_\zeta W(\xi,\eta,\zeta)
&=&
B(\xi)\zeta.
\end{eqnarray*}
We suppose that $\Gamma_D$ is a nonempty
Lipschitz portion of $\partial \Omega$ and that
a bounded trace operator 
$\gamma_0:~H^1(\Omega)\rightarrow L^2(\Gamma_D)$ is defined and 
\[
V:=\{v\in H^1(\Omega));~\gamma_0(v)=0\},
\]
is not empty.
Then the minimization problem \ref{mini2} corresponds to
the following linear elliptic boundary value problem.
\[
-\DIV (B(x)\nabla u)+b(x)u=f(x)
~~~~~
\mbox{in}~\Omega(\varphi),
\]
\[
u=g~~~\mbox{on}~\Gamma_D,
~~~~~
(B\nabla u)\cdot\nu=0~~~\mbox{on}~\partial\Omega\setminus\Gamma_D,
\]
where $\nu$ is a unit normal vector on $\partial\Omega$.
From the Poincar\'{e} inequality and the Lax-Milgram theorem,
there exists a unique global minimizer $u(\varphi)\in V(g,\varphi)$
for $\varphi\in\cO(\Omega)$.
}
\end{Ex}
\begin{Ex}\label{lapex}
{\rm
We consider the following potential energy
which corresponds to the Poisson equation.
\begin{equation}\label{linearprob2}
W(\xi,\eta,\zeta)=\half |\zeta |^2-f(\xi)\eta,
\end{equation}
\[
E(v,\Omega)=\int_\Omega
\left(
\half
|\nabla v|^2-f(x)v
\right)dx,
\]
where $k\in\N\cup\{0\}$ and $f\in W^{k,2}(\Omega_0,\R)$.
This is a special case of Example~\ref{linex}.
We remark that
\[
\nabla_\xi W(\xi,\eta,\zeta)
=
-\nabla f(\xi)\eta,~~~~~
W_\eta (\xi,\eta,\zeta)
=
-f(\xi),~~~~~
\nabla_\zeta W(\xi,\eta,\zeta)
=
\zeta .
\]
Under the same boundary condition, the minimization problem \ref{mini2} 
corresponds to
the following boundary value problem of the Poisson equation.
\[
-\Delta u=f(x)
~~~~~
\mbox{in}~\Omega(\varphi),
\]
\[
u=g~~~\mbox{on}~\Gamma_D,
~~~~~
\frac{\partial u}{\partial \nu}=0~~~\mbox{on}~\partial\Omega\setminus\Gamma_D.
\]
In the case of $n=2$, this represents the anti-plane
displacement model of the isotropic linear elasticity.
}
\end{Ex}
\begin{Lem}\label{lemE}
For the potential energy of Example~\ref{linex}, 
$\E\in C^k(H^1(\Omega)\times \cO(\Omega))$ and 
{\rm (\ref{dEformula})}
holds.
\end{Lem}
\noindent
{\it Proof.}
In the case of Example~\ref{linex}, $\E$ becomes
\begin{eqnarray*}
\E(v,\varphi)
&=&
\int_\Omega
\left\{
\half\left(A(\varphi)\nabla v\right)^\srT 
B(\varphi(x))
\left(A(\varphi)\nabla v\right)
+\half
b(\varphi(x))v^2
-f(\varphi(x))v
\right\}
\kappa(\varphi)dx\\
&=&
\int_\Omega
\left\{
\half\Psi_1(v,\varphi)(x)+\half\Psi_2(v,\varphi)(x)
-\Psi_3(v,\varphi)(x)
\right\}
\kappa(\varphi)(x)dx
\end{eqnarray*}
where
\[
\Psi_1(v,\varphi):=
\left(A(\varphi)\nabla v\right)^\srT 
(B\circ\varphi)
\left(A(\varphi)\nabla v\right),
~~~
\Psi_2(v,\varphi):=
(b\circ\varphi) v^2,
~~~
\Psi_3(v,\varphi):=
(f\circ \varphi)v.
\]
Under the assumptions, from Theorem~\ref{kappaA}, we obtain
\[
[(v,\varphi)\mapsto A(\varphi)\nabla v]
\in C^\infty(H^1(\Omega)\times \cO(\Omega),L^2(\Omega)^n),
\]
\[
[\varphi\mapsto B\circ\varphi]
\in C^k(\cO(\Omega),L^\infty(\Omega,\R^{n\times n})),
\]
\[
[(w,B)\mapsto w^\srT Bw]\in
C^\infty(L^2(\Omega)^n\times L^\infty(\Omega,\R^{n\times n}),
L^1(\Omega)).
\]
From these regularities,
it follows that
$\Psi_1\in C^k(H^1(\Omega)\times \cO(\Omega),L^1(\Omega))$.
Similarly, from the following regularities
\[
[\varphi\mapsto b\circ\varphi]
\in C^k(\cO(\Omega),L^\infty(\Omega)),
~~~
[(b,v)\mapsto bv^2]\in
C^\infty(L^\infty (\Omega)\times H^1(\Omega),L^1(\Omega)),
\]
\[
[\varphi\mapsto f\circ\varphi]
\in C^k(\cO(\Omega),L^2(\Omega)),
~~~
[(f,v)\mapsto fv]\in
C^\infty(L^2 (\Omega)\times H^1(\Omega),L^1(\Omega)),
\]
$\Psi_2$ and $\Psi_3$
also belong to $C^k(H^1(\Omega)\times \cO(\Omega),L^1(\Omega))$.
Since $\kappa\in C^\infty (\cO(\Omega),L^\infty(\Omega))$
from Theorem~\ref{kappaA}, we conclude that
$\E\in C^k(H^1(\Omega)\times \cO(\Omega))$.
\qed 

We note that if $\Omega$ satisfies the cone property,
the regularity conditions (\ref{bfcond1})
can be weakened by using the Sobolev imbedding theorem 
(see \cite{Ada75} etc.):
$H^1(\Omega)$ is continuously imbedded in $L^p(\Omega)$,
where $p\in [1,\infty)$ if $n=2$ and $p\in [1,2n/(n-2)]$
if $p\geq 3$. The relaxed conditions are
\begin{equation}\label{bfcond2}
~^\exists q>1,~~~
b\in W^{k,q}(\Omega_0,\R),
~~~~~
f\in W^{k,q}(\Omega_0,\R),~
\mbox{if $n=2$,}
\end{equation}
\begin{equation}\label{bfcond3}
b\in W^{k,\frac{n}{2}}(\Omega_0,\R),
~~~~~
f\in W^{k,\frac{2n}{n+2}}(\Omega_0,\R),~
\mbox{if $n\geq 3$.}
\end{equation}
More precisely, we have the following proposition. 
\begin{Prop}
We suppose that $\Omega$ satisfies the cone property.
For the potential energy of Example~\ref{linex}
with the condition {\rm (\ref{bfcond2})} or 
{\rm (\ref{bfcond3})} in stead of {\rm (\ref{bfcond1})},
$\E\in C^k(H^1(\Omega)\times \cO(\Omega))$ and 
{\rm (\ref{dEformula})}
holds.
\end{Prop}
\noindent
{\it Proof.}
We suppose $n=2$ and the condition (\ref{bfcond2}).
For the $q>1$ in (\ref{bfcond2}) we define
$q^*:=(1-1/q)^{-1}$. 
From
\[
[\varphi\mapsto b\circ\varphi]
\in C^k(\cO(\Omega),L^q(\Omega)),
~~~
[v\mapsto v]\in C^\infty(H^1(\Omega),L^{2q^*}(\Omega)),
\]
\[
[(b,v)\mapsto bv^2]\in
C^\infty(L^q (\Omega)\times L^{2q^*}(\Omega),L^1(\Omega)),
\]
$\Psi_2\in C^k(H^1(\Omega)\times \cO(\Omega),L^1(\Omega))$
follows. In the same way, from
\[
[\varphi\mapsto f\circ\varphi]
\in C^k(\cO(\Omega),L^q(\Omega)),
~~~
[v\mapsto v]\in C^\infty(H^1(\Omega),L^{q^*}(\Omega)),
\]
\[
[(f,v)\mapsto fv]\in
C^\infty(L^q (\Omega)\times L^{q^*}(\Omega),L^1(\Omega)),
\]
$\Psi_3\in C^k(H^1(\Omega)\times \cO(\Omega),L^1(\Omega))$
follows.

Next we suppose $n\geq 3$ and the condition (\ref{bfcond3}).
Put $p:=n/(n-2)$,\, $p^*:=(1-1/p)^{-1}=n/2$,
and $(2p)^*:=(1-1/(2p))^{-1}=2n/(n+2)$. Then we have
\[
[v\mapsto v]\in C^\infty(H^1(\Omega),L^{2p}(\Omega)),
\]
\[
[\varphi\mapsto b\circ\varphi]
\in C^k(\cO(\Omega),L^{p^*}(\Omega)),
~~~
[(b,v)\mapsto bv^2]\in
C^\infty(L^{p^*} (\Omega)\times L^{2p}(\Omega),L^1(\Omega)),
\]
\[
[\varphi\mapsto f\circ\varphi]
\in C^k(\cO(\Omega),L^{(2p)^*}(\Omega)),
~~~
[(f,v)\mapsto fv]\in
C^\infty(L^{(2p)^*}(\Omega)\times L^{2p}(\Omega),L^1(\Omega)).
\]
Hence, from these regularities, we obtain
$\Psi_2,~\Psi_3\in C^k(H^1(\Omega)\times \cO(\Omega),L^1(\Omega))$.
\qed 

We state our results under 
the condition (\ref{bfcond1})
hereafter. But we remark that the following results are valid 
even under the condition (\ref{bfcond2}) or 
(\ref{bfcond3}) with the cone property.
We obtain the following theorem from Lemma~\ref{lemE} and
Theorem~\ref{ckcondition}.
\begin{Th}\label{main1}
For the potential energy of Example~\ref{linex}, 
$E_*\in C^k(\cO(\Omega))$ holds, and if $k\geq 1$ we have
\[
E_*'(\varphi_0)[\mu]=
\int_\Omega
\left(
\nabla_\xi W(u)\cdot \mu 
-(\nabla_\zeta W(u))^\srT (\nabla\mu^\srT )\nabla u
+W(u)\DIV \mu
\right)dx,
\]
for $\mu\in W^{1,\infty}(\Omega_0,\R^n)$,
where $u$ is the global minimizer of $E(\cdot,\Omega)$ in $V(g)$.
In particular, in the case of Example~\ref{lapex}, we have
\[
E_*'(\varphi_0)[\mu]=
\int_\Omega
\left\{
-(\nabla f\cdot \mu)u
-(\nabla^\srT u)(\nabla\mu^\srT )(\nabla u)
+\left(
\half |\nabla u|^2-fu
\right)\DIV \mu
\right\}dx,
\]
for $\mu\in W^{1,\infty}(\Omega_0,\R^n)$.
\end{Th}

We consider the case of Example~\ref{linex}. The global minimizer $u$ 
belongs to $H^1(\Omega)$, but not to $H^2(\Omega)$ in general.
The following boundary integral (\ref{limf}) is called J-integral
(\cite{Oht85} etc.).
\begin{Th}
Under the assumptions in Example~\ref{linex} with $k=1$,
if there exists a sequence of subdomains $\{\Omega_l\}_l$ 
in which the Gauss-Green formula
holds and
\[
\Omega_1\subset \Omega_2\subset \cdots \Omega
\mbox{ with }\bigcup_{l=1}^\infty\Omega_l=\Omega,
\]
and if the global minimizer $u$ belongs to $H^2(\Omega_l)$
for each $l\in\N$, 
then we have
\begin{equation}\label{limf}
E_*'(\varphi_0)[\mu]=
\lim_{l\rightarrow \infty}
\int_{\partial \Omega_l}
\left\{
W(u)~\mu\cdot\nu
-
\left(\nabla_\zeta W(u)\cdot\nu \right)
\left(\nabla u\cdot\mu \right)\right\}
d\cH^{n-1}_x,
\end{equation}
for $\mu\in W^{1,\infty}(\Omega_0,\R^n)$,
where $\nu$ denotes the outward unit normal of $\partial \Omega_l$
and $\cH^{n-1}$ denotes the $n-1$ dimensional Hausdorff measure.
\end{Th}
\noindent
{\it Proof.}
Under the conditions, we have 
\[
-\DIV \left[\nabla_\zeta W(u)\right]
+W_\eta(u)
=0~~~~~\mbox{in}~L^2(\Omega_l),
\]
which is equivalent to
\[
-\DIV (B\nabla u)+bu=f
~~~~~\mbox{in}~L^2(\Omega_l).
\]
Let us denote the Hessian matrix of $u$ by $\nabla^2 u$.
Since
\begin{eqnarray*}
\int_{\Omega_l} W(u(x))\DIV \mu(x) dx
&=&
\int_{\partial \Omega_l}W(u(x))\mu(x)\cdot\nu d\cH^{n-1}_x
-\int_{\Omega_l}\nabla_x [W(u(x))]\cdot\mu(x)dx,
\end{eqnarray*}
\[
\nabla_x [W(u(x))]
=
\nabla_\xi W(u(x))+W_\eta (u(x))\nabla u(x)
+[\nabla^2 u(x)]\nabla_\zeta W(u(x)),
\]
\begin{eqnarray*}
\int_{\Omega_l}W_\eta (u(x))\nabla u(x)\cdot\mu(x)dx
&=&
\int_{\Omega_l} \DIV \left[\nabla_\zeta W(u)\right]
\nabla u(x)\cdot\mu(x)dx\\
&=&
\int_{\partial \Omega_l}
\left(\nabla_\zeta W(u)\cdot\nu \right)
\left(\nabla u(x)\cdot\mu(x)\right)
d\cH^{n-1}_x\\
&&
-\int_{\Omega_l} \nabla_\zeta^\srT  W(u)
\nabla\left(\nabla u(x)\cdot\mu(x)\right)dx\\
&=&
\int_{\partial \Omega_l}
\left(\nabla_\zeta W(u)\cdot\nu \right)
\left(\nabla u(x)\cdot\mu(x)\right)
d\cH^{n-1}_x\\
&&
-\int_{\Omega_l} \nabla_\zeta^\srT  W(u)
\left\{(\nabla^2 u(x))\mu(x)+(\nabla\mu^\srT )\nabla u\right\}dx,
\end{eqnarray*}
\begin{eqnarray*}
\int_{\Omega_l}\nabla_x [W(u(x))]\cdot\mu(x)dx
&=&
\int_{\Omega_l}\left\{
\nabla_\xi W(u(x))
+[\nabla^2 u(x)]\nabla_\zeta W(u(x)) \right\}\cdot\mu(x)dx\\
&&\int_{\partial \Omega_l}
\left(\nabla_\zeta W(u)\cdot\nu \right)
\left(\nabla u(x)\cdot\mu(x)\right)
d\cH^{n-1}_x\\
&&
-\int_{\Omega_l} \nabla_\zeta^\srT  W(u)
\left\{(\nabla^2 u(x))\mu(x)+(\nabla\mu^\srT )\nabla u\right\}dx\\
&=&
\int_{\partial \Omega_l}
\left(\nabla_\zeta W(u)\cdot\nu \right)
\left(\nabla u(x)\cdot\mu(x)\right)
d\cH^{n-1}_x\\
&&
+\int_{\Omega_l}\left\{
\nabla_\xi W(u)
 \cdot\mu dx
-\nabla_\zeta^\srT  W(u)
(\nabla\mu^\srT )\nabla u\right\}dx, \\
 & &
\end{eqnarray*}
we obtain, for $\mu\in W^{1,\infty}(\Omega_0,\R^n)$,
\begin{eqnarray*}
E_*'(\varphi_0)[\mu]
&=&
\int_\Omega
\left(
\nabla_\xi W(u)\cdot \mu 
-(\nabla_\zeta W(u))^\srT (\nabla\mu^\srT )\nabla u
+W(u)\DIV \mu
\right)dx\\
&=&
\int_{\Omega\setminus \Omega_l}
\left(
\nabla_\xi W(u)\cdot \mu 
-(\nabla_\zeta W(u))^\srT (\nabla\mu^\srT )\nabla u
+W(u)\DIV\mu
\right)dx\\
&&
+\int_{\Omega_l}
\left(
\nabla_\xi W(u)\cdot \mu 
-(\nabla_\zeta W(u))^\srT (\nabla\mu^\srT )\nabla u
\right)dx\\
&&+\int_{\partial \Omega_l}W(u)~\mu\cdot\nu ~d\cH^{n-1}_x
-\int_{\Omega_l}\nabla_x [W(u)]\cdot\mu dx\\
&=&
\int_{\Omega\setminus \Omega_l}
\left(
\nabla_\xi W(u)\cdot \mu 
-(\nabla_\zeta W(u))^\srT (\nabla\mu^\srT )\nabla u
+W(u)\DIV\mu
\right)dx\\
&&
+\int_{\partial \Omega_l}
\left\{
W(u)~\mu\cdot\nu
-
\left(\nabla_\zeta W(u)\cdot\nu \right)
\left(\nabla u\cdot\mu \right)\right\}
d\cH^{n-1}_x .
\end{eqnarray*}
The first term tends to $0$ if $l\rightarrow \infty$.
Hence, we have the formula (\ref{limf}).
\qed

In particular, for the case of Example~\ref{lapex}, we have
\begin{eqnarray*}
E_*'(\varphi_0)[\mu]
&=&
\lim_{l\rightarrow \infty}
\int_{\partial \Omega_l}
\left\{
\left(
\half |\nabla u|^2-fu
\right)~\mu \cdot\nu
-
\left(\nabla u\cdot\nu \right)
\left(\nabla u\cdot\mu\right)\right\}
d\cH^{n-1}_x .
\end{eqnarray*}

\begin{Th}
Let $\Omega$ be a bounded domain with $C^2$-boundary.
Under the condition of Example~\ref{lapex}, we assume that
$\Gamma_D=\partial \Omega$ with $g\equiv 0$.
Then the following formula holds.
\[
E_*'(\varphi_0)[\mu]
=
-\half\int_{\partial \Omega}|\nabla u|^2\mu\cdot\nu ~d\cH^{n-1}_x .
\]
\end{Th}
\noindent
{\it Proof.}
From the regularity theorem of elliptic boundary value problems,
we have $u\in H^2(\Omega)$.
So, we can choose $\Omega_l=\Omega$.
We remark that 
\[
\nabla u\cdot \mu =(\nabla u\cdot\nu)(\mu\cdot \nu),
~~~~~
|\nabla u\cdot\nu|=|\nabla u|
~~~~~\mbox{on}~\partial \Omega,
\]
since $u$ and its tangential derivatives vanish on the boundary. 
Hence we have
\begin{eqnarray*}
E_*'(\varphi_0)[\mu]
&=&
\int_{\partial \Omega}
\left\{
\left(
\half |\nabla u|^2-fu
\right)~\mu \cdot\nu
-
\left(\nabla u\cdot\nu \right)
\left(\nabla u\cdot\mu\right)\right\}
d\cH^{n-1}_x\\
&=&
\int_{\partial \Omega}
\left\{
\left(
\half |\nabla u|^2
\right)~\mu \cdot\nu
-
\left(
|\nabla u|^2
\right)~\mu \cdot\nu\right\}
d\cH^{n-1}_x\\
&=&
-\half \int_{\partial \Omega}
|\nabla u|^2
~\mu \cdot\nu~
d\cH^{n-1}_x.
\end{eqnarray*}
\qed

This theorem can be generalized as follows.
\begin{Th}
Under the condition of Example~\ref{lapex}, we assume that
$\mu\in W^{1,\infty}(\Omega_0,\R^n)$ satisfies 
\[
\supp (\mu)\cap \supp(~g~|_{\ov{\Gamma_D}})=\emptyset,
~~~~~
\supp (\mu)\cap \ov{\Gamma_D}\cap 
\ov{\partial\Omega\setminus\Gamma_D}=\emptyset,
\]
\[
\mbox{$~^\exists \cO$ : an open set of $\R^n$ s.t.
$\cO\supset \supp(\mu)$ and $\partial\Omega\cap\cO$ is $C^2$-class.
}
\]
Then the following formula holds.
\[
E_*'(\varphi_0)[\mu]
=
-\half\int_{\Gamma_D}|\nabla u|^2\mu\cdot\nu ~d\cH^{n-1}_x 
-\int_{\partial\Omega\setminus\Gamma_D}f\, u~\mu\cdot\nu ~d\cH^{n-1}_x.
\]
\end{Th}

\end{document}